\documentclass[a4paper,12pt]{article}
\usepackage{amsfonts}
\usepackage{amsmath}
\usepackage{amssymb}
\usepackage{amsthm}
\usepackage{url}
\title{Zeta functions of supersingular curves \\of genus 2}
\author{Daniel Maisner$^{*}$ and Enric Nart\thanks{acknowledge
financial support from the project BFM-2003-06092 of MCYT}}
\date{}
\newtheorem{teor}{Theorem}[section]
\newtheorem{cor}[teor]{Corollary}

\newtheorem{rem}[teor]{Remark}
\newtheorem{lem}[teor]{Lemma}
\newtheorem{prop}[teor]{Proposition}

\newcommand{\N}{\mathbb N}

\newcommand{\Z}{\mathbb Z}

\newcommand{\C}{\mathbb C}
\def\op{\operatorname}

\def\as{\op{AS}}

\def\be{\bigskip}

\def\c#1#2{{\mathcal C}_{(#1,#2)}}

\def\diso{\lower.4ex\hbox{$\downarrow$}\raise.4ex\hbox{\mbox{\scriptsize $\wr$}}}
\def\e{\medskip}

\def\eps{\epsilon}
\def\exp{\op{exp}}

\def\ff#1{\mathbb F_{q^{#1}}}

\def\fq{\mathbb F_q}
\def\ft{\mathbb F_2}
\def\f4{\mathbb F_4}

\def\gen#1{\big\langle\, {#1} \,\big\rangle}

\def\imp{\ \Longrightarrow\ }

\def\k{\op{Ker}}
\def\kb{\bar{k}}

\def\la{\lambda}
\def\lg{l\raise.6ex\hbox to.2em{\hss.\hss}l}
\def\lra{\longrightarrow}

\def\md#1{\ \mbox{\rm(mod }{#1})}

\def\ni{\noindent}

\def\orb{\hbox to  .3em{$\backslash$}\backslash}

\def\sg{\sigma}

\def\sii{\,\Longleftrightarrow\,}

\def\sn{\op{sgn}}

\def\tr{\op{Tr}}

\def\wb{\overline{W}}

\newcounter{cs}
\stepcounter{cs}
\newcommand{\casos}{\begin{itemize}}
\newcommand{\fcasos}{\end{itemize}\setcounter{cs}{1}}

\newfont{\tit}{cmr12 scaled \magstep3}

\begin{document}
\maketitle

\section*{Introduction}
This paper was motivated by the problem of determining what
isogeny classes of abelian surfaces over a finite field $k$
contain jacobians. In \cite{mn} we performed a numerical
exploration of this problem, that led to several conjectures. We
present in this paper a complete answer for supersingular surfaces
in charac\-teristic 2 (section \ref{isg}). We deal with this
problem in a direct way by computing explicitly the zeta function
of all supersingular curves of genus two (section \ref{prfxd}).
Our procedure is constructive, so that we are able to exhibit
curves with prescribed zeta function and to count the number of
curves, up to $k$-isomorphism, leading to the same zeta function.

We base our work on the ideas of van der Geer and van der Vlugt
\cite{vv1},\,\cite{vv2}, who expressed the number of points of a
supersingular curve of genus two in terms of certain invariants.
In section \ref{vw} we compute explicitly these invariants in
terms of the coefficients of a defining equation and in section
\ref{zeta} we compute the number of points of the curve over the
quadratic extension in terms of objects defined over $k$.

\section{Supersingular curves of genus 2 in characteristic 2}\label{as}
In this section we review the results of van der Geer-van der
Vlugt and we fix some notations. Let $k=\fq$ be a finite field of
even characteristic, with $q=2^m$. We recall some basic facts
concerning the Artin-Schreier operator:

$$ \as\colon k\lra k,\quad\as(x)=x+x^2.
$$

This is an $\ft$-linear operator with kernel $\ft$. The image
$\as(k)$ is an $\ft$-subspace of $k$ of codimension one; hence,
$|\as(k)|=q/2$ and $k/\as(k)\simeq \ft$.

We shall denote simply by $\tr$ or $\tr_k$ the absolute trace
$\tr_{k/\ft}$. For any $x\in k$ we have $\tr(x)=\tr(x^2)$, because
$x^2$ is a galois conjugate of $x$ over the prime field $\ft$.
Therefore, $\as(k)=\k(\tr)$.

For any $a\in k$, the polynomial $x^2+x+a\in k[x]$ is separable.
Its roots are in $k$ if and only if $a\in \as(k)$. Hence,

\begin{lem}\label{qeq}
A quadratic polynomial $f(x)=x^2+ax+b\in k[x]$ is separable iff
$a\ne0$; in this case, $f(x)$ is irreducible iff
$b/a^2\not\in\as(k)$.
\end{lem}

Throughout the paper we shall denote by $\mu_n\subseteq \kb$ the
group of the $n$-th roots of 1 and we let $\eps\in\mu_3$ be a
fixed root of the polynomial $x^2+x+1$. Also, we denote
$k_n:=\ff{n}$. Note that $k\subseteq \as(k_2)$, since
$\tr_{k_2/k}(k)=0$. Clearly,
$$\renewcommand\arraystretch{1.1}
\begin{array}{ccccl}
1\in\as(k)&\iff& \f4\subseteq k&\iff& m \mbox{ even},\\
\mu_3\subseteq k&\iff& (k^*)^3\varsubsetneq k^*&\iff& m\mbox{
even},\\ \mu_3\subseteq\as(k)&\iff&
\mathbb{F}_{16}\subseteq k&\iff& m\equiv 0\md4,\\
\mu_5\subseteq k&\iff& (k^*)^5\varsubsetneq k^*&\iff&
m\equiv0\md4.
\end{array}\renewcommand\arraystretch{1.} $$

\e

Every projective smooth curve of genus 2, defined over $k$ and
supersingular (i.e. with supersingular jacobian) admits an affine
model of the type:
\begin{equation}\label{eqss}
C\colon\quad y^2+y=ax^5+bx^3+cx+d,\quad a\in k^*,\,b,\,c\in
k,\,d\in k/\as(k),
\end{equation}
which has only one point at infinity. We can think that the term
$d$ takes only two values, $d=0$ or $d=d_0$, with $d_0\in
k-\as(k)$ fixed. To apply the {\it hyperelliptic twist} to the
curve $C$ consists in adding $d_0$ to the defining equation. If we
denote by $C^{\tau}$ the twisted curve, we have
\begin{equation}\label{tw} |C(\fq)|+|C^{\tau}(\fq)|=2q+2.
\end{equation}
The curves $C$ and $C^{\tau}$ are isomorphic over the quadratic
extension of $k$ through the mapping $(x,y)\mapsto (x,y+u)$, where
$u\in k_2$ satisfies $u+u^2=d_0$.

Throughout the paper we abuse of language and identify the curve
$C$ given by (\ref{eqss}) with the family $(a,b,c,d)$ of the four
parameters involved in the defining equation.

\begin{rem}
The mappings $(x,y)\mapsto (x,y+cx)$, $(x,y)\mapsto
(x,y+cx+c^2x^2)$ set respective $k$-isomorphisms between the curve
(\ref{eqss}) and the curves $$ y^2+y=ax^5+bx^3+c^2x^2+d,\qquad
y^2+y=ax^5+c^4x^4+bx^3+d,
$$which are the models used respectively in \cite{vv1} and
\cite{cnp}. We beg the reader to pay attention to this change of
models when we quote results of these two papers.
\end{rem}

By (\ref{tw}), in order to study the number of points of these
curves we can assume $d=0$. Consider the linear polynomial
$R(x)=ax^4+bx^2+c^2x$. Since $\tr(cx)=\tr(c^2x^2)$, the function:
$$
Q\colon k\lra\ft,\quad x\mapsto \tr(ax^5+bx^3+cx)=\tr(xR(x)),
$$
is a quadratic form associated to the simplectic form:
$$
k\times k\lra \ft,\quad (x,y)\mapsto
\gen{x,y}_R=\tr\left(xR(y)+yR(x)\right),
$$
since clearly, \begin{equation} \label{qq}
Q(x+y)=Q(x)+Q(y)+\gen{x,y}_R,\quad\forall x,y\in k.
\end{equation}

The number of zeros of $Q$ determines the number of points of $C$:
$$
|C(\fq)|=1+2|Q^{-1}(0)|.
$$

The radical of the simplectic form $\gen{\ ,\ }_R$ coincides with
the set of roots in $k$ of the $\ft$-linear and separable
polynomial (independent of $c$):
$$
E_{ab}(x):=a^4x^{16}+b^4x^8+b^2x^2+ax.
$$
Let $\wb=\k(E_{ab})$ denote the subspace of $\kb$ formed by the 16
roots of this polynomial. We denote
$$
W:=\op{rad}\gen{\ ,\ }_R=\wb \cap k,\quad w:=\dim_{\ft}(W),\quad
0\le w\le 4.
$$

From (\ref{qq}) we deduce:
$$
Q(x+y)=Q(x)+Q(y), \quad \forall x\in k,\, y\in W.
$$
In particular, $Q$ defines a linear form, $Q\colon W\lra \ft$. The
space $V:=\k(Q_{|W})$ controls the behaviour of $Q$ on the classes
$x+W$; for all $x\in k$, $y\in W$:
$$
Q(x+y)=Q(x) \iff y\in V.
$$
This subspace $V$ of $W$ has codimension 0 or 1. If
$V\varsubsetneq W$, in each class $x+W$ the quadratic form $Q$
vanishes on half of the elements; therefore $|Q^{-1}(0)|=q/2$ and
$|C(\fq)|=1+q$.

If $V=W$, the quadratic form $Q$ is constant on each class $x+W$.
Hence, it factorizes through a quadratic form, which we still
denote by $Q$,
$$
Q\colon k/W\lra \ft,
$$
associated to the non-degenerate simplectic form induced by
$\gen{\ ,\ }_R$ on $k/W$. In particular, the dimension of $k/W$ is
even, so that $m$, $w$ have the same parity. Moreover, if
$m-w=2n$, the number of zeros of $Q$ can take only two values:
$2^{n-1}(2^n+1)$ or $2^{n-1}(2^n-1)$. Thus,
$$
|C(\fq)|=1+2\left(2^w(2^{n-1}(2^n\pm1))\right)=1+q\pm \sqrt{2^wq}.
$$

Summarizing,\e

\renewcommand\arraystretch{1.2}
\ni{\bf Theorem (van der Geer-van der Vlugt). }
$$
\begin{array}{l}
 V\varsubsetneq W \imp |C(\fq)|=1+q,\\
 V=W\imp |C(\fq)|=1+q\pm \sqrt{2^wq}.
\end{array}
$$\renewcommand\arraystretch{1.}

There are, thus, three invariants that determine the number of
points of $C$: the dimension $w$ of the space $W$, the codimension
0 or 1 of the subspace $V\subseteq W$ and the sign ``$+$" or
``$-$" telling the parity, even or odd, of the quadratic form $Q$,
in the case $V=W$. Actually, the last two invariants can be
unified using the following terminology:
$$
\sn(Q):=\left\{\begin{array}{ll} 0,\quad &\mbox{ if
}V\varsubsetneq W,\\+/-,\quad &\mbox{ the parity of
$Q_{\,|(k/W)}$, if }V=W.
\end{array}\right.
$$

We end this section of preliminaries recalling the conditions that
are ne\-cessary and sufficient for two models (\ref{eqss}) to give
$k$-isomorphic curves. In general (cf. \cite[Lemma 2.3]{vv2} or
\cite[Proposition 10]{cnp}), the supersingular curves given
respectively by $(a,b,c,d)$, $(a',b',c',d')$ are $k$-isomorphic if
and only if there exist $\la\in k^*$, $\nu\in k$ such that:
\begin{equation}\label{modiso}
(a',b',c',d')=\left(\la^5a,\la^3b,\la\left(c+\root
4\of{E_{ab}(\nu)}\right), a\nu^5+b\nu^3+c\nu+d\right),
\end{equation}
the equality $d'=a\nu^5+b\nu^3+c\nu+d$ understood in $k/\as(k)$.

There are $4q-2+[8]_{4\mid m}$ $k$-isomorphism classes of
supersingular curves \cite[Theorem 2]{cnp}, where $+[8]_{4\mid m}$
means ``add $8$ if $4\mid m$".

The two cases $b=0$, $b\ne0$, give disjoint families of
isomorphism classes of supersingular curves. The curves with $b=0$
are all isomorphic to the curve $y^2+y=x^5$ over $\kb$ and they
have 160 automorphisms; the curves with $b\ne 0$ are
$\kb$-isomorphic to curves of the type $y^2+y=a(x^5+x^3)$ and have
32 automorphisms.

Note that if $b\ne0$, we can achieve $\la^5a=\la^3b$  by taking
$\la=\sqrt{b/a}$. Hence, in this case we can always assume that
$a=b$.

As a consequence of (\ref{modiso}) we see that any curve with
$|C(\fq)|=q+1$ is isomorphic to its own hyperelliptic twist. In
fact, since $V\varsubsetneq W$, there exists $\nu\in W$ with
$Q(\nu)\ne0$, that is, $a\nu^5+b\nu^3+c\nu\not\in\as(k)$; hence,
the curve $(a,b,c,0)$ is isomorphic to the curve $(a,b,c,d_0)$.

\section{Computation of the invariants $W$, $V$}\label{vw}
In this section we compute explicitly the subspaces $W$, $V$  in
terms of the parameters $a,\,b,\,c$ of the defining equation of
the curve.

\subsection*{Computation of $W$}
The polynomial $E_{ab}$ factorizes in $k[x]$ \cite[Theorem
3.4]{vv1}:
\begin{multline*}
E_{ab}(x)=a^4x^{16}+b^4x^8+b^2x^2+ax=\\
x(a^2x^5+b^2x+a)(a^2x^{10}+b^2x^6+ax^5+1)=xP(x)(1+x^5P(x)),
\end{multline*}
 with $P(x):=P_{ab}(x):=a^2x^5+b^2x+a$. Hence, we have $v^5P(v)=1,0$ respectively
 for 10,6 elements $v\in\wb$.

\begin{lem}\label{ptbase}
\begin{enumerate}
\item Any family of 4 roots of $P(x)$ is a basis of $\wb$.
\item The 10 elements $v\in\wb$ such that $vP(v)=1$ can be
expressed in a unique way as the sum of two roots of $P(x)$.
\end{enumerate}
\end{lem}

\begin{proof}
Let $z_1,\,z_2,\,z_3,\,z_4,\,z_5\in\kb$ be the  roots of $P(x)$
and let us check that $z_1,\,z_2,\,z_3,\,z_4$ are linearly
independent. They are all non-zero and different, hence, the sum
of any two of them cannot vanish. Since $z_1+z_2+z_3+z_4+z_5=0$,
the sum of any three or four of them cannot vanish either.

The 10 elements $z\in\wb$, such that $vP(v)=1$ are the sum of two
or three of the elements of the basis. In any case, they are the
sum of two roots of $P(x)$, uniquely determined.
\end{proof}

\begin{lem}\label{av5}
Let $v=z+z'$ be a root in $k$ of $v^5P(v)=1$, with $z,z'$ roots of
$P(x)$. Then,
$$
\begin{array}{l}
(av^5)^{-1}\in\as(k)\imp z,z'\in k,\\
(av^5)^{-1}\not\in\as(k)\imp z,z'\in k_2-k \mbox{ and they are
conjugate over }k.
\end{array}
$$
\end{lem}

\begin{proof}
Let us impose that $v+z$ is a root of $P(x)$:
\begin{multline*}
0=a^2(v+z)^5+b^2(v+z)+a=a^2v^5+a^2v^4z+a^2vz^4+a^2z^5+b^2v+b^2z+a=\\
=a^2v^5+a^2v^4z+a^2vz^4+b^2v=v(a^2v^4+a^2v^3z+a^2z^4+b^2).
\end{multline*}
We deduce that $a^2v^4+a^2v^3z+a^2z^4+b^2=0$. If we multiply by
$z$ and apply $a^2z^5+b^2z=a$, we get
$$a^2v^3z^2+a^2v^4z+a=0,$$ or equivalently,
$z^2+vz+(av^3)^{-1}=0$. Since $v\ne 0$, this equation in $z$ is
separable and the two roots are $z,z'$. By Lemma \ref{qeq}, the
roots belong to $k$ iff $(av^5)^{-1}\in\as(k)$.
\end{proof}

We are ready to see that the factorization of $P(x)$ as a product
of irreducible polynomials determines $w$. We shall write
$P(x)=(n_1)(n_2)\cdots(n_t)$ to indicate that $P(x)$ factorizes in
$k[x]$ as the product of $t$ irreducible polynomials of degrees
$n_1,\,n_2,\dots,n_t$.

\begin{prop}\label{wp} Let $P(x)=a^2x^5+b^2x+a$, with $a\in k^*$, $b\in k$. Then,
\begin{enumerate}
\item $w=0\iff $ $P(x)$ is irreducible. \item $w=1\iff $
$P(x)=(1)(4)$ or $P(x)=(2)(3)$. \item $w=2\iff $ $P(x)=(1)(1)(3)$
or $P(x)=(1)(2)(2)$. \item $w=3\iff $ $P(x)=(1)(1)(1)(2)$.
\item $w=4\iff $ $P(x)=(1)(1)(1)(1)(1)$.
\end{enumerate}
\end{prop}

\begin{proof}
If $P(x)$=(1)(1)(1)(1)(1), we have $W=\wb$ by Lemma \ref{ptbase}
and $w=4$.

If $P(x)$=(1)(1)(1)(2), we have $W\varsubsetneq\wb$ and $W$
contains the 3 roots of $P(x)$ in $k$, which are linearly
independent by Lemma \ref{ptbase}. Hence, $w=3$.

Suppose that $P(x)$=(1)(1)(3) and let $z,z'$ be the roots of
$P(x)$ in $k$. By Lemma \ref{av5}, $z+z'$ is the only root of
$1+x^5P(x)$ that belongs to $k$. Hence, $W$ is the subspace
generated by $z,z'$ and $w=2$. Suppose now that $P(x)$=(1)(2)(2).
By Lemma \ref{av5}, the two traces of the quadratic factors of
$P(x)$ are the only roots of $1+x^5P(x)$ that belong to $k$. Thus,
 $W$ has 4 elements and $w=2$.

Similarly, by Lemma \ref{av5} we have $w=1$ if $P(x)$=(1)(4) or
$P(x)$=(2)(3), and we have $w=0$ if $P(x)$ is irreducible.

Since we have considered all possible factorizations of $P(x)$,
the converse implications hold too.
\end{proof}

We proceed now to find explicit criteria to determine the
factorization type of $P(x)$ in terms of $a,\,b$. We start with an
 auxiliary result.

\begin{lem}\label{deesc}
Let $e\in k^*$. The polynomial $f(x)=x^4+x^3+x^2+x+e$ has the
following decomposition in $k[x]$ as a product of irreducible
factors:\renewcommand\arraystretch{1.2}
$$
\begin{array}{lcll}
e\not\in (k^*)^3&&&\imp\ f(x)=(1)(3),\\
e=\la^3,\,&\la\in k-\as(k),\,&m\mbox{ odd}
&\imp\ f(x)\mbox{ is irreducible},\\
e=\la^3,\,&\la\in \as(k),\,&m\mbox{ odd}&\imp\ f(x)=(1)(1)(2),\\
e=\la^3,\,&\la\mu_3\nsubseteq \as(k),\,&m\mbox{ even}&\imp\ f(x)=(2)(2),\\
e=\la^3,\,&\la\mu_3\subseteq \as(k),\,&m\mbox{ even}&\imp\
f(x)=(1)(1)(1)(1).
\end{array}
$$
\end{lem}
\renewcommand\arraystretch{1.}

\begin{proof}
We check first that $e=\la^3$, $\la\in\as(k)$, are necessary and
sufficient conditions in order that $f(x)$ decomposes in $k[x]$ as
the product of two polynomials of degree 2, not necessarily
irreducible. In fact, assume that we have such a decomposition:
\begin{equation}\label{desr}
 x^4+x^3+x^2+x+e=(x^2+ux+s)(x^2+(u+1)x+t);
\end{equation}
this amounts to:
$$
s+t=1+u+u^2,\ u(s+t)+s=1,\ st=e.
$$
From the first and second equations we deduce $s=(u+1)^3$,
$t=u^3$, so that $e=(u+u^2)^3$. Conversely, if $e=\la^3$ and
$\la=u+u^2$, with $\la,u\in k$, we get the decomposition above by
taking $s=(u+1)^3$, $t=u^3$.

By Lemma \ref{qeq}, the quadratic factor $x^2+ux+(u+1)^3$ is
irreducible iff $(u+1)^3/u^2$ does not belong to $\as(k)$. Since:
$$
(u+1)^3/u^2=u+1+u^{-1}+u^{-2},
$$
this condition is equivalent to $u+1\not\in\as(k)$. Similarly, the
quadratic factor $x^2+(u+1)x+u^3$ is irreducible iff
$u\not\in\as(k)$.

We start now the proof of the lemma. Assume that $e\not\in
(k^*)^3$. Then, $x^3+e$ is irreducible in $k_2[x]$. Therefore,
$f(x)$=(1)(3), since this is the only factorization for which
$f(x)$ is not the product of two polynomials of degree 2 over
$k_2$.

If $m$ is odd, then $e=\la^3$ for a unique $\la\in k$. If
$\la\not\in \as(k)$, then $f(x)$ does not factorize as the product
of two polynomials of degree 2 in $k[x]$, but it admits such a
factorization over $k_2[x]$; hence, $f(x)$ is irreducible. On the
other hand, if $\la=u+u^2$, with $u\in k$, we have a factorization
(\ref{desr}). Now, since $m$ is odd, we have $1\not\in\as(k)$ and
necessarily $f(x)$=(1)(1)(2), since exactly one of the two
conditions, $u+1\not\in\as(k)$, $u\not\in\as(k)$, is satisfied.

Suppose now that $e\in(k^*)^3$ and $m$ is even. If $\la^3=e$, with
$\la\in k$, the elements $\la\eps$ and $\la\eps^2$ are cubic roots
of $e$ too. Since their sum is zero, either all three belong to
$\as(k)$, or only one of them. This corresponds to $f(x)$ having
three different decompositions (\ref{desr}) or only one, that is,
to $f(x)$=(1)(1)(1)(1) or $f(x)=$(2)(2).
\end{proof}

In order to study the decomposition of $P_{ab}(x)$ we can assume
that $b=0$ or $b=a$, as remarked at the end of section \ref{as}.

\renewcommand\arraystretch{1.4}
\begin{prop}\label{beso}
Let $a\in k^*$ and $P_{a0}(x)=a^2(x^5+a^{-1})$. Then,
$$
P_{a0}(x)=\left\{
\begin{array}{ll}
(1)(4),&\mbox{ if $m$ is odd},\\
(1)(2)(2),&\mbox{ if }m\equiv 2\md4,\\
(1)(1)(1)(1)(1), &\mbox{ if }m\equiv0\md4,\ a\in (k^*)^5,\\
\mbox{irreducible}, &\mbox{ if }m\equiv0\md4,\ a\not\in (k^*)^5.
\end{array}
\right.
$$
\end{prop}
\renewcommand\arraystretch{1.}

\begin{proof}
Suppose first $a\not\in (k^*)^5$, or equivalently, that
$P_{a0}(x)$ has no roots in $k$. We have necessarily
$m\equiv0\md4$ and $\mu_5\subseteq k$. Thus, if we adjoint to $k$
any root of $P_{a0}(x)$, this polynomial will split completely in
the larger field. Thus, $P_{a0}(x)$ cannot be (2)(3) and  it must
be irreducible.

Suppose now $a\in(k^*)^5$ and let $z\in k$ satisfy $z^5=a^{-1}$.
We have $$ x^5+a^{-1}=(x+z)(x^4+zx^3+z^2x^2+z^3x+z^4),
$$
and the quartic factor has the same factorization type as the
polynomial $x^4+x^3+x^2+x+1$, that has been studied in Lemma
\ref{deesc}.
\end{proof}

\renewcommand\arraystretch{1.2}
\begin{prop}\label{bneso}
Let $a\in k^*$ and $P_{aa}(x)=a^2(x^5+x+a^{-1})$. Suppose that
$P_{aa}(x)$ has no roots in $k$. Then,

$$
P_{aa}(x)=\left\{
\begin{array}{ll}
(2)(3),&\mbox{ if $m$ odd},\\
\mbox{irreducible}, &\mbox{ if $m$ even}.
\end{array}
\right.
$$

Suppose that $P_{aa}(x)$ has a root $z\in k$. Then, for
$e=1+z^{-4}$, we have:

$$
P_{aa}(x)=\left\{ \begin{array}{lcll}
(1)(1)(3),& \mbox{ if }e\not\in (k^*)^3,&&\\
(1)(4)& \mbox{ if }e=\la^3,\,&\la\in k-\as(k),\,&m\mbox{ odd},\\
(1)(1)(1)(2)&\mbox{ if }e=\la^3,\,&\la\in \as(k),\,&m\mbox{ odd},\\
(1)(2)(2)&\mbox{ if }e=\la^3,\,&\la\mu_3\nsubseteq \as(k),\,&m\mbox{ even} ,\\
(1)(1)(1)(1)(1)&\mbox{ if }e=\la^3,\,&\la\mu_3\subseteq
\as(k),\,&m\mbox{ even}.
\end{array}\right.
$$
\end{prop}
\renewcommand\arraystretch{1.}

\begin{proof}
If the polynomial has no roots in $k$, the assertion is
consequence of Proposition \ref{wp} and the fact that $m$ and $w$
have the same parity.

If the polynomial has some root $z\in k$, then,
$$
x^5+x+a^{-1}=(x+z)(x^4+zx^3+z^2x^2+z^3x+z^4+1),
$$
and the quartic factor has the same decomposition type as the
polynomial $x^4+x^3+x^2+x+(1+z^{-4})$, that has been obtained in
Lemma \ref{deesc}
\end{proof}

This result allows us to count the number of times that appears
each decomposition of $P_{aa}(x)$, when $a$ varies. This
computation is crucial to find in section \ref{prfxd} explicit
formulas for the number of curves with prescribed zeta function.

\begin{cor}\label{nbra}
The two following tables give the number of values of $a\in k^*$
leading to each of the
 possible factorizations of $P_{aa}(x)$,
respectively in the cases $m$ odd and $m$ even:\be

\renewcommand\arraystretch{1.4}
\centerline{
\begin{tabular}{ccc}
\hline $(2)(3)$&$(1)(1)(1)(2)$&$(1)(4)$\\
\hline$(q+1)/3$&$(q-2)/6$&$(q/2)-1$
\end{tabular}
}\be \centerline{
\begin{tabular}{cccc}
\hline $(1)(1)(3)$&$(1)(2)(2)$&$(1)(1)(1)(1)(1)$&irreducible\\
\hline$(q-1)/3$&$(q/4)-[1]_{4\nmid m}$&$((q-4)/60)-\left[\frac
15\right]_{4\mid m}$&$\frac 25(q+1-[2]_{4\mid m})$
\end{tabular}
}\be
\renewcommand\arraystretch{1.}
\end{cor}

\begin{proof}
Suppose that $m$ is odd. By Proposition \ref{bneso}, the values of
$a\in k^*$ leading to $P_{aa}(x)=$(1)(4) are parameterized by
elements $\la\in k-\as(k)$, $\la\ne 1$, via
\begin{equation}\label{lza}
1+z^{-4}=\la^3,\quad a=(z^5+z)^{-1}.
\end{equation}
These two relations set $\la$ in 1-1 correspondence with $z$
(since $(k^*)^3=k^*$) and $z$ in 1-1 correspondence with $a$
(since $z^5+z=a^{-1}$ has only one root). We get $(q/2)-1$ values
of $a$.

Similarly, the values of $a\in k^*$ leading to
$P_{aa}(x)=$(1)(1)(1)(2) are parameterized by choosing
$\la\in\as(k)$, $\la\ne 0$ and taking $z,\,a$ as before. The
relation between $\la$ and $z$ is still 1-1, but now there are
three different values of $z$ linked to the same $a$. We get $1/3$
of the values computed above.

All other values of $a\in k^*$ lead to $P_{aa}(x)=$(2)(3).\e

Suppose now $m$ even. There are $2(q-1)/3$ values of $z\in k$
satisfying $1+z^{-4}\not\in (k^*)^3$, and each two of these values
give the same $a=(z^5+z)^{-1}$. We have thus $(q-1)/3$ values of
$a$ with $P_{aa}(x)=$(1)(1)(3).

In order to ensure the factorization $P_{aa}(x)=$(1)(2)(2), we
take $\la\not\in \as(k)\cup \f4$ and  consider $z,\,a$ determined
by (\ref{lza}). The fact that $\la\not\in\f4$ guarantees that
$z\ne0,1$ and  $a\ne 0$. The number of different values of $\la$
with these properties is:
\begin{equation}\label{quan1} \frac
q2\,-2, \ \mbox{ if }4\nmid m,\quad \frac q2, \ \mbox{ if }4\mid
m.
\end{equation}  Since $\la+\la\eps+\la\eps^2=0$,
in the couple $\la\eps$, $\la\eps^2$ exactly one element belongs
to $\as(k)$. The element not belonging to $\as(k)$ and $\la$ give
the same value of $a$. Hence, the number of values of $a$ is half
of the quantities given in (\ref{quan1}).

For $P_{aa}(x)$ to split completely, we have to take $\la\in k$
such that $\la,\,\la\eps\in \as(k)$; this will ensure that
$\la\mu_3\subseteq \as(k)$. By the non-degeneracy of the pairing
$\tr(xy)$, there are $q/4$ values of $\la\in k$ with this
property: $\la\in\gen{1,\eps}^{\perp
}=\f4^{\perp}$. Also, we need
$\la\not\in\f4$  in order that $z\ne 0,1$. Since
$\f4^{\perp}\cap\f4=\{0\}$ if $4\nmid m$ and $\f4^{\perp}\supseteq
\f4$ if $4\mid m$, the number of values of $\la$ is,
\begin{equation}\label{quan2}
\frac q4\,-1 \ \mbox{ if }4\nmid m,\quad \frac q4\,-4, \ \mbox{ if
}4\mid m.
\end{equation} Every 3 values of $\la$ give the same $z$
and every 5 values of $z$ give the same $a$. Hence, the number of
different values of $a$ is obtained dividing by 15 the numbers
given in (\ref{quan2}).

All other values of $a\in k^*$ lead to $P_{aa}(x)$ irreducible.
\end{proof}

\subsection*{Computation of $V$}
We can use Lemmas \ref{ptbase} and \ref{av5} to reinterpret the
linear form $Q_{|W}$ in a way that provides an explicit
computation of $\op{codim}(V,W)$. Let us start with some remarks
on linear forms over $k$. For any $c\in k$, let us denote by $L_c$
the linear form,
$$
L_c\colon k\lra \ft,\quad x\mapsto \tr_k(cx).
$$
The non-degeneracy of the pairing $\tr(xy)$ allows us to consider
a linear isomorphism:
$$
L\colon k\lra \op{Hom}(k,\ft),\quad c\mapsto L_c
$$
In particular, for any subspace $W\subseteq k$ of dimension $w$,
the linear mapping,
$$
L\colon k\lra \op{Hom}(W,\ft),\quad c\mapsto {L_c}_{\,|W},
$$
is onto and each linear form over $W$ has $q/2^w$ preimages.

\renewcommand\arraystretch{1.2}
\begin{prop}\label{qsenzill}
Let $(a,b,c,d)$ be parameters defining a supersingular curve
(\ref{eqss}). Let $\ell,\ell_c\colon W\lra\ft$ be the linear forms
determined by:
$$
\begin{array}{ll}
\ell(z)=\tr(1),\quad &  \mbox{ if }P(z)=0,\\\ell(v)=0,\quad &
\mbox{ if }v=z+z',\mbox{ with $z,z'$ roots of $P(x)$ in }k,\\
\ell(v)=1,\quad & \mbox{ if }v=z+z',\mbox{ with $z,z'$ roots of
$P(x)$ in }k_2-k,
\end{array}
$$and  $\ell_c=L_{c+b^2a^{-1}}$ restricted to $W$.

Then, $Q_{|W}=\ell_c+\ell$. In particular, $V=W$ iff
$\ell_c=\ell$.
\end{prop}
\renewcommand\arraystretch{1.}

\begin{proof}
Suppose that $P(z)=0$, with $z\in k$. We have $a^2z^5+b^2z+a=0$.
If we multiply by $z^5$, we get $az^5+a^2z^{10}=b^2z^6$, so that
$b^2z^6\in\as(k)$ and, in consequence, $bz^3\in \as(k)$. We can
now compute:
$$
Q(z)=\tr(az^5+bz^3+cz)=\tr(b^2a^{-1}z+1+cz)=\ell_c(z)+\tr(1).
$$

If $v^5P(v)=1$, we have $a^2v^{10}+b^2v^6+av^5+1=0$, so that
$b^2v^6\equiv 1\md{\as(k)}$ and, in consequence,  $bv^3\equiv
1\md{\as(k)}$. On the other hand,
$$
av^5=b^2a^{-1}v+1+(av^5)^{-1},
$$
so that,
$$
av^5+bv^3\equiv b^2a^{-1}v+(av^5)^{-1} \md{\as(k)}.
$$
By Lemma \ref{av5}, $\tr((av^5)^{-1})=0,1$ according to $z,z'$
belonging to $k$ or to $k_2-k$. This ends the proof.
\end{proof}

Note that $\ell$ depends on $a,b$ and $\ell_c$ depends on $a,b,c$.
Thus, for $a,b$ fixed the invariant $\op{codim(V,W)}$ is
determined by the linear form $\ell_c$, or equivalently, by the
linear form ${L_c}_{|W}$. Let us check now that this linear form
determines the $k$-isomorphism class of the curve too, up to
hyperelliptic twist.

\begin{lem}\label{nbrclas}
Let $C=(a,b,c,0)$ and suppose that $b\ne 0$ or $4\nmid m$. Then,
for any $c'\in k$, the following conditions are equivalent:
\begin{enumerate}
\item[i)] $\left\{\begin{array}{ll} C'=(a,b,c',0) \mbox{ is
$k$-isomorphic to $C$},&\mbox{ if }V\varsubsetneq
W,\\C'=(a,b,c',0) \mbox{ is $k$-isomorphic either to $C$ or to
$C^{\tau}$},&\mbox{ if }V=W.\end{array}\right.$ \item[ii)] $c'\in
c+\root 4\of{E_{ab}(k)}$, \item[iii)]
${L_c}_{|W}={L_{c'}}_{|W}$,\item[iv)] $\ell_c=\ell_{c'}$,
\end{enumerate}
\end{lem}

\begin{proof}
Conditions (i) and (ii) are equivalent by (\ref{modiso}), since
our hypothesis on $b$ and/or $m$ imply that $\lambda=1$ in
(\ref{modiso}).

In order to check that (ii) and (iii) are equivalent, let us show
that $E_{ab}(k)=\left(W^4\right)^{\perp}$, where the orthogonal is
taken with respect to the isomorphism of $k$ with its dual
obtained from the perfect pairing $\tr(xy)$. In fact, for an
arbitrary $\nu\in k$ we have:
\begin{equation}\label{amnt}
0=\gen{z,\nu}_R=\tr(z(a\nu^4+b\nu^2)+\nu(az^4+bz^2)),\quad\forall
z\in W.
\end{equation}
Clearly,
$$
z(a\nu^4+b\nu^2)\equiv z^4(a^4\nu^{16}+b^4\nu^8)\md{\as(k)},\quad
\nu bz^2\equiv \nu^2 b^2z^4\md{\as(k)},
$$
so that (\ref{amnt}) is equivalent to:
\begin{equation}\label{avall}
\tr(z^4E_{ab}(\nu))=0,\ \forall z\in W.
\end{equation}Hence,
$E_{ab}(k)\subseteq (W^4)^{\perp}$ and, having the same dimension,
they coincide. Therefore, condition (ii) is equivalent to $c'\in
c+W^{\perp}$, which is is equivalent to (iii) by the definition of
$L_c$.

Finally, it is obvious that (iii) and (iv) are equivalent.
\end{proof}

\section{Computation of the zeta function}\label{zeta}
Let $C$ be a smooth projective curve of genus 2, defined over $k$.
The zeta function of $C$ is a formal series in one indeterminate,
which can be expressed as a rational function:
\begin{equation}\label{a1a2}
\op{Z}(C/\fq,t)=\exp\left(\sum_{n\ge1}
N_n\frac{t^n}n\right)=\frac{1+a_1t+a_2t^2+qa_1t^3+q^2t^4}{(1-t)(1-qt)},
\end{equation}
where $N_n:=|C(\ff{n})|$ and $a_1,a_2\in\Z$. From this identity
one deduces immediately that:
$$ N_1=q+1+a_1,\quad N_2=q^2+1+2a_2-a_1^2.
$$ Thus, the zeta function is determined by the couple $(N_1,N_2)$.

Let $J_C$ be the jacobian variety of $C$. The polynomial
$t^4+a_1t^3+a_2t^2+qa_1t+q^2$ is the characteristic polynomial of
the Frobenius endomorphism of the abelian surface $J_C$. This
polynomial determines the $k$-isogeny class of $J_C$; thus, two
curves have the same zeta function if and only if their jacobian
varieties are $k$-isogenous.

Let $C$ be a supersingular curve defined by (\ref{eqss}), with
parameters $(a,b,c,0)$. Denote by $w$, $V$, $W$, $Q$, $\ell$,
$\ell_c$ the objects associated to $C_{|k}$ in sections
\ref{as},\,\ref{vw} and by $\tilde{w}$, $\tilde{V}$, $\tilde{W}$,
$\tilde{Q}$, $\tilde{\ell}$, $\tilde{\ell}_c$ the corresponding
objects associated to the curve $C_{|k_2}$.

In this section we compute $N_2=|C(\ff2)|$ in terms of $a,b,c$.
The idea is to apply the results of the last section and take
advantage of the fact that the curve is defined over $k$ to avoid
any computation in $k_2$. More precisely, we shall see that the
linear form $\ell_c$ on $W$ contains already sufficient
information to determine $N_2$.

To begin with we recall some observations on linear forms over
$k_2/k$. For any $c\in k$ we denote by $\tilde{L}_c$ the linear
mapping:
$$
\tilde{L}_c\colon k_2\lra \ft,\quad x\mapsto
\tilde{L}_c(x)=\tr_{k_2}(cx).
$$
As before, we can consider the isomorphism,
$$
\tilde{L}\colon k\lra \op{Hom}(k_2/k,\ft), \quad c\mapsto
\tilde{L}_c.
$$

\renewcommand\arraystretch{.6}
\begin{lem}\label{diagtr}
Think of $\tr_{k_2/k}$ as a linear isomorphism between $k_2/k$ and
$k$ and denote by $(\tr_{k_2/k})^*$ its dual isomorphism. We have
a commutative diagram of linear isomorphisms:
$$
\begin{array}{rcl}
\op{Hom}(k,\ft)&\stackrel{(\tr_{k_2/k})^*}{\lra}&\op{Hom}(k_2/k,\ft)\\&&\\
\qquad L\nwarrow&&\nearrow\tilde{L}\qquad\\
&\lower1ex \hbox{$k$}&
\end{array}
$$
\end{lem}
\renewcommand\arraystretch{1.}

\begin{proof}
By the transitivity of the trace, for any $x\in k_2$ we have
$$
\tilde{L}_c(x)=\tr_{k_2}(cx)=\tr_k(\tr_{k_2/k}(cx))=\tr_k(c\tr_{k_2/k}(x))=
L_c(\tr_{k_2/k}(x)).
$$
\end{proof}

The invariant $\tilde{w}$ is completely determined by the
factorization of $P(x)$ in $k[x]$, which was obtained in
Propositions \ref{beso},\, \ref{bneso}. The invariant
$\op{codim(\tilde{V},\tilde{W})}$ can be determined as follows:

\begin{prop}\label{v=w}
$\tilde{V}=\tilde{W}$ iff the following two conditions are
satisfied:
\begin{enumerate}
\item $\ell_c(v)=0$, for $v=z+z'\in W$, with $z,z'$
roots of $P(x)$ in $k_2-k$. \item $\ell_c(z)=1$, if $P(x)=(1)(4)$
over $k[x]$ and $z$ is its root in $k$.
\end{enumerate}
\end{prop}

\begin{proof}
Assume first that $\tilde{V}=\tilde{W}$. By Proposition
\ref{qsenzill}, we have $\tilde{\ell}_c=\tilde{\ell}$. If
$v=z+z'\in k$, with $z,z'$ roots of $P(x)$ in $k_2-k$, we have
$\tilde{\ell}(z)=0$ by definition; hence, by Lemma \ref{diagtr},
$\ell_c(v)=\tilde{\ell}_c(z)=\tilde{\ell}(z)=0$. If $P(x)=(1)(4)$,
then it factorizes over $k_2$ as:
$P(x)=(x+z)(x^2+ux+t)(x^2+u'x+t')$, with $u,u'\in k_2-k$ galois
conjugate and $z+u+u'=0$. By definition, $\tilde{\ell}(u)=1$;
hence, by Lemma \ref{diagtr},
$\ell_c(z)=\tilde{\ell}_c(u)=\tilde{\ell}(u)=1$.

Assume now that conditions (1), (2) are satisfied and let us check
that $\tilde{\ell}_c=\tilde{\ell}$. For any root $z$ of $P(x)$ in
$k_2$ we have $\tilde{\ell}(z)=0$; if $z\in k$ we have directly
$\tilde{\ell}_c(z)=0$, whereas for $z\in k_2-k$ with galois
conjugate $z'$ we have $\tilde{\ell}_c(z)=\ell_c(z+z')=0$ by
condition (1). In particular, if $v=z+z'$, with $z,z'$ roots of
$P(x)$ in $k_2$, we have $\tilde{\ell}(v)=0$ and
$\tilde{\ell}_c(v)=0$ too. Finally, $\tilde{\ell}(v)=1$ if
$v=\omega+\omega'$, with $\omega,\omega'$ roots of $P(x)$ in
$k_4-k_2$; in this case necessarily $P(x)=(1)(4)$ over $k[x]$,
$z:=\tr_{k_2/k}(v)$ is the only root of $P(x)$ in $k$ and
$\tilde{\ell}_c(v)=\ell_c(z)=1$ by condition (2).
\end{proof}

We address now to the computation of the sign of $\tilde{Q}$ when
$\tilde{V}=\tilde{W}$. The crucial observation is that over $k_2$
we have:
$$\gen{k+\tilde{W},k+\tilde{W}}_R=0,\quad \tilde{Q}(k+\tilde{W})=0,$$
since $k\subseteq \as(k_2)$ and $\tilde{Q}(\tilde{W})=0$ by
assumption. This will allow us to control the behavior of
$\tilde{Q}$ on the classes of elements of $k_2$ modulo
$k+\tilde{W}$.

The simplectic form  $\gen{\ ,\ }_R$ is non-degenerate over
$k_2/\tilde{W}$; hence,
\begin{multline*}
\dim\left((k+\tilde{W})/\tilde{W}\right)^{\perp}=\dim
k_2/\tilde{W}-\dim(k+\tilde{W})/\tilde{W}=\\
=(2m-\tilde{w})-(m-w)=m+w-\tilde{w}.
\end{multline*}
Let $k+\tilde{W}\subseteq U\subseteq k_2$ be the subspace such
that $U/\tilde{W} =\left((k+\tilde{W})/\tilde{W}\right)^{\perp}$.
We know that $\dim U=m+w$. Clearly, $\gen{\ ,\ }_R$ induces a
non-degenerate simplectic form:
\begin{equation}\label{u}
U/(k+\tilde{W})\,\times\,U/(k+\tilde{W})\,\lra\,\ft,\quad
(x,y)\mapsto \gen{x,y}_R,
\end{equation}
on the space $U/(k+\tilde{W})$, of dimension $2w-\tilde{w}$. Let
$n:=w-(\tilde{w}/2)$.

For arbitrary $x\in k_2$, $y\in k+\tilde{W}$ we have:
\begin{equation}\label{qu}
\tilde{Q}(x+y)=\tilde{Q}(x)+\tilde{Q}(y)+\gen{x,y}_R=\tilde{Q}(x)+\gen{x,y}_R.
\end{equation}
For fixed $x$, the linear mapping,
$$
k+\tilde{W}\lra \ft,\quad y\mapsto \gen{x,y}_R
$$
vanishes only for $x\in U$. Thus, for $x\in U$, $\tilde{Q}$ is
constant in the class $x+(k+\tilde{W})$ and it determines a
quadratic form $\tilde{Q}\colon U/(k+\tilde{W})\lra \ft$
associated to the simplectic form (\ref{u}). The number of zeros
of $\tilde{Q}$ will be $2^{n-1}(2^n\pm 1)$ and altogether there
are
$$2^{n-1}(2^n\pm 1)2^{m+\tilde{w}-w}=2^{m+w-1}\pm 2^{m+(\tilde{w}/2)-1}
$$ zeros of $\tilde{Q}$ in $U$.

 Moreover, for $x\not\in U$,
 $\gen{x,-}_R$ does not vanish on $k+\tilde{W}$ and, by (\ref{qu}), $\tilde{Q}$ takes the values 0,1 the same
 number of times, $2^{m+\tilde{w}-w-1}$, in the class $x+(k+\tilde{W})$. There
 are  $2^{\dim k_2/(k+\tilde{W})}-2^{\dim U/(k+\tilde{W})}=2^{m+w-\tilde{w}}-2^{2w-\tilde{w}}$ such classes and we count
$$\left(2^{m+w-\tilde{w}}-2^{2w-\tilde{w}}\right)2^{m+\tilde{w}-w-1}=2^{2m-1}-2^{m+w-1}
$$ zeros of $\tilde{Q}$ in $k_2-U$.

Therefore, the number of points of $C$ over $k_2$ is:
$$
|C(\ff2)|=1+2\left(2^{2m-1}\pm
2^{m+(\tilde{w}/2)-1}\right)=1+q^2\pm \sqrt{2^{\tilde{w}}q^2}.
$$
We have thus proved,

\begin{prop}\label{sgnq}
If $\tilde{V}=\tilde{W}$ the sign of $\tilde{Q}$ as a quadratic
form over $k_2/\tilde{W}$ coincides with the sign of $\tilde{Q}$
as a quadratic form over $U/(k+\tilde{W})$.
\end{prop}

In order to determine this latter sign, we find an explicit
description of $U$ and we express the action of $\tilde{Q}$ on $U$
in terms of the action of $Q$ on $W$.

\begin{prop}\label{sgnu}
We have $U=\tr_{k_2/k}^{-1}(W)$. Moreover, for any $u\in U$, with
relative trace $z=\tr_{k_2/k}(u)\in W$, we have,
$$\tilde{Q}(u)=Q(z) \sii P(z)=0 \mbox{ or }z=0.$$
\end{prop}

\begin{proof}
For any $u\in k_2$ with $\tr_{k_2/k}(u)=z$, we have $u\in U$ iff:
\begin{equation}\label{amunt}
0=\gen{u,\la}_R=\tr_{k_2}\left(\la(au^4+bu^2)+u(a\la^4+b\la^2)\right),\quad\forall
\la\in k.
\end{equation}
By the same argument used to show that (\ref{amnt}) was equivalent
to (\ref{avall}), we see that (\ref{amunt}) is equivalent to
\begin{multline*}
\tr_{k_2}(\la^4E_{ab}(u))=0,\ \forall \la\in k \sii
\tr_k(\la^4E_{ab}(z))=0,\ \forall \la\in k\iff\\\sii \tr_k(\la
E_{ab}(z))=0,\ \forall \la\in k \sii E_{ab}(z)=0,
\end{multline*}
the last equivalence by the non-degeneracy of the pairing
$\tr(xy)$. This proves the first assertion.

Now, let $u\in U$. The galois conjugate of $u$ is $u^{\sg}=u+z$.
Hence,
\begin{multline*}
\tr_{k_2/k}(au^5+bu^3+cu)=au^5+bu^3+cu+a(u+z)^5+b(u+z^3)+c(u+z)=\\
=au^4z+auz^4+az^5+bu^2z+buz^2+bz^3+cz=\\=az^5+bz^3+cz+
uR(z)+zR(u),
\end{multline*}
so that,
\begin{multline*}
\tilde{Q}(u)=\tr_{k_2}(au^5+bu^3+cu)=\tr_k(az^5+bz^3+cz+
uR(z)+zR(u))=\\=Q(z)+\tr_k(uR(z)+zR(u)).
\end{multline*}
We have to check when $uR(z)+zR(u)\in\as(k)$. Let us express
$u=zv$, with $v\in k_2$ an element of relative trace 1: $v^2+v=r$,
$r\in k-\as(k)$. Note that $v^4=v^2+r^2=v+r+r^2$. We have,
\begin{multline*}
uR(z)+zR(u)=zv(az^4+bz^2)+z(az^4v^4+bz^2v^2)=\\
=v(az^5+bz^3)+v^4az^5+v^2bz^3
=(r+r^2)az^5+rbz^3\equiv\\\equiv(r+r^2)az^5+r^2b^2z^6=
raz^5+r^2(a^2z^{10}+z^5P(z))\equiv
\\\equiv r^2z^5P(z)\md{\as(k)}.
\end{multline*}
Hence, if $z^5P(z)=0$ we have $uR(z)+zR(u)\in\as(k)$ and if
$z^5P(z)=1$ we have $uR(z)+zR(u)\equiv
 r^2 \not\equiv 0\md{\as(k)}$, since $r\not\in \as(k)$.
\end{proof}

\begin{teor}\label{taulasn}
The possible signs of $\tilde{Q}$ are given in the following
table:\be

\renewcommand\arraystretch{1.2}
\centerline{
\begin{tabular}{|c|c|c|c|c|}
\hline $w$&$\tilde{w}$&$P(x)$&$\dim U/(k+\tilde{W})$&$\sn(\tilde{Q})$\\
\hline $0$&$0$&irreducible&$0$&$+$\\
\hline $1$&$2$&$(1)(4)$ or $(2)(3)$&$0$&$0/+$\\
\hline $2$&$2$&$(1)(1)(3)$&$2$&$+/-$\\
\hline $2$&$4$&$(1)(2)(2)$&$0$&$0/+$\\
\hline $3$&$4$&$(1)(1)(1)(2)$&$2$&$0/+/-$\\
\hline $4$&$4$&$(1)(1)(1)(1)(1)$&$4$&$+/-$\\
\hline
\end{tabular}
}\renewcommand\arraystretch{1.} \be

Moreover, let $Z\subseteq W$ be the subset of all roots of $P(x)$
in $k$. For $P(x)=(1)(1)(3)$ or $(1)(1)(1)(2)$ we have,
$$
\sn(\tilde{Q})=``-" \sii \ell_c(z)\ne\tr(1), \ \forall z\in Z,
$$whereas for $P(x)=(1)(1)(1)(1)(1)$ we have,
$$\sn(\tilde{Q})=``+" \sii \ell_c(z)=0, \ \mbox{for exactly 3 of the 5 roots }z\in Z.
$$
\end{teor}

\begin{proof}
The content of the table is an immediate consequence of
Propositions \ref{wp}, \ref{v=w} and \ref{sgnq}. The other
assertions on $\sn{\tilde{Q}}$ are consequence of Propositions
\ref{sgnq} and \ref{sgnu}. For instance, in the cases when $\dim
U/(k+\tilde{W})=2$ the quadratic form $\tilde{Q}$ has either 1 or
3 zeros on this space; the minus sign correspond to the case
$\tilde{Q}(u)=1$ for all $u\in U/(k+\tilde{W})$, $u\ne0$, and by
Proposition \ref{sgnu} this is equivalent to $Q(z)=1$ for all
$z\in Z$, which is equivalent to $\ell_c(z)=\tr(1)+1$ for all
$z\in Z$ by Proposition \ref{qsenzill}.

We leave the case $w=\tilde{w}=4$ to the reader.
\end{proof}

We have obtained an explicit computation of the zeta function of
any supersingular curve, except for the sign of $Q$ when $V=W$.
From the computational point of view, once you know $\pm a_1$ and
$a_2$, the sign of $a_1$ is easy to determine by computing
iterates of a random divisor in the jacobian. We can consider a
deterministic algorithm too by evaluating $Q$ on a simplectic
basis of $k/W$ with respect to $\gen{\-,\-}_R$.

\section{Zeta functions of supersingular curves of genus 2}\label{prfxd}
In this section we compute all possible zeta
functions arising from supersingular curves of genus 2 and we find
formulas for the number of $k$-isomorphism classes of curves that
have the same zeta function. We proceed in a constructive way, by
applying the results of the previous sections to all supersingular
curves; hence, our results can be used to exhibit curves with
prescribed zeta function.

For any couple of integers $(a_1,a_2)$ we shall denote by
$\c{a_1}{a_2}$ the set of $k$-isomorphism classes of smooth
projective curves of genus 2 defined over $k$, whose zeta function
is given by (\ref{a1a2}), or equivalently, whose number of points
$N_1$, $N_2$ over $k$ and $k_2$ satisfy $$ N_1=q+1+a_1,\quad
N_2=q^2+1+2a_2-a_1^2.
$$The hyperelliptic twist sets a bijection between $\c{a_1}{a_2}$
and $\c{-a_1}{a_2}$, which is the identity if $a_1=0$ by the
remark at the end of section \ref{as}.

We work with supersingular curves given by equation (\ref{eqss}),
depending on four parameters $(a,b,c,d)$ with $b=0$ or $b=a$. We
keep the notations $w$, $W$, $V$, $Q$, $\ell$, $\tilde{w}$,
$\tilde{W}$, $\tilde{V}$, $\tilde{Q}$ introduced in the last
section. We remind that
$$
\ell_c={L_c}_{\,|W}, \ \mbox{ if }b=0,\quad
\ell_c={L_{c+a}}_{\,|W}, \ \mbox{ if }b=a.
$$
We deal first with the case $m$ odd.

\subsection{$\mathbf{P_{ab}(x)=(1)(4)}$}
By Proposition \ref{wp} we have $w=1$, $\tilde{w}=2$ in this case.
If $z\in k$ is the only root of $P_{ab}(x)$ in $k$ we have
$W=\{0,z\}$.

Let us study first the case $b=0$. Since $(k^*)^5=k^*$, we can
assume $a=1$ by (\ref{modiso}), so that $z=1$. By Proposition
\ref{beso}, $P_{10}(x)=(1)(4)$. By Propositions \ref{qsenzill},
\ref{v=w} and Theorem \ref{taulasn}, for any $c\in k$ we have
\begin{equation}\label{snqq}\renewcommand\arraystretch{1.2}
\begin{array}{l}
\ell_c(z)=0\imp \sn(Q)=\sn(\tilde{Q})=0,\\ \ell_c(z)\ne0 \imp
\sn(Q)=\pm, \ \sn(\tilde{Q})=+.
\end{array}\renewcommand\arraystretch{1.}
\end{equation}By Lemma \ref{nbrclas}, the
different values of $(c,d)$ lead to three $k$-isomorphism classes
represented by $(1,0,0,0)$ and a couple of twisted curves,
$(1,0,1,0)$, $(1,0,1,1)$. The first one has
$(N_1,N_2)=(q+1,q^2+1)$ and the other two $(N_1,N_2)=(q+1\pm
\sqrt{2q},q^2+1+2q)$. Thus, we get one curve in each of the sets
$\c{0}{0}$, $\c{\sqrt{2q}}{2q}$, $\c{-\sqrt{2q}}{2q}$.

Let us study now the case $a=b\ne0$. By Corollary \ref{nbra},
there are $(q/2)-1$ values of $a$ leading to $P_{aa}(x)=(1)(4)$
and by (\ref{modiso}) they correspond to different $k$-isomorphism
classes. Let us fix one of these values of $a$. As before,
(\ref{snqq}) holds and the different values of $(c,d)$ provide
three $k$-isomorphism classes, represented by $(a,a,0,0)$,
$(a,a,c,0)$, $(a,a,c,1)$, where $\ell_c(z)\ne0$, and they are
distributed into the same three zeta functions.

We have altogether a contribution of $q/2$ $k$-isomorphism classes
in each of the sets $\c{0}{0}$, $\c{\sqrt{2q}}{2q}$,
$\c{-\sqrt{2q}}{2q}$.

\subsection{$\mathbf{P_{ab}(x)=(2)(3)}$}
By Proposition \ref{wp} we have $w=1$, $\tilde{w}=2$ in this case.
If $x^2+vx+t$ is the quadratic irreducible factor of $P_{ab}(x)$
we have $W=\{0,v\}$.

By Proposition \ref{beso}, we have necessarily $b\ne0$ and we can
assume $a=b$. By Corollary \ref{nbra} we have $(q+1)/3$ values of
$a$ leading to this factorization of $P_{aa}(x)$. We fix one of
these values of $a$. By Propositions \ref{qsenzill}, \ref{v=w} and
Theorem \ref{taulasn}, for any $c\in k$ we have
$$\renewcommand\arraystretch{1.2}
\begin{array}{l}
\ell_c(v)=0\imp \sn(Q)=0,\ \sn(\tilde{Q})=+,\\ \ell_c(v)\ne0 \imp
\sn(Q)=\pm, \ \sn(\tilde{Q})=0
\end{array}\renewcommand\arraystretch{1.}
$$
By Lemma \ref{nbrclas}, the different values of $(c,d)$ lead to
three $k$-isomorphism classes represented by $(a,a,a,0)$,
$(a,a,c,0)$, $(a,a,c,1)$, where $\ell_c(v)\ne0$. The first one has
$(N_1,N_2)=(q+1,q^2+1+2q)$ and the other two $(N_1,N_2)=(q+1\pm
\sqrt{2q},q^2+1)$. Thus, we get $(q+1)/3$ curves in each of the
sets $\c{0}{q}$, $\c{\sqrt{2q}}{q}$, $\c{-\sqrt{2q}}{q}$.

\subsection{$\mathbf{P_{ab}(x)=(1)(1)(1)(2)}$}
By Proposition \ref{wp} we have $w=3$, $\tilde{w}=4$ in this case.
We have $W=\gen{z_1,z_2,z_3}$, where $z_1,z_2,z_3$ are the roots
of $P_{ab}(x)$ in $k$. The quadratic irreducible factor of
$P_{ab}(x)$ is $x^2+vx+t$, with $v=z_1+z_2+z_3$.

By Proposition \ref{beso}, we have necessarily $b\ne0$ and we
assume $a=b$. By Corollary \ref{nbra} we have $(q-2)/6$ values of
$a$ leading to this factorization of $P_{aa}(x)$. For any such
fixed value of $a$, the linear form $\ell\colon W\lra \ft$
introduced in Proposition \ref{qsenzill} is determined by
$\ell(z_1)=\ell(z_2)=\ell(z_3)=1$.

For any $c\in k$, let $N$ be the number of $z_i$ such that
$\ell_c(z_i)=0$. Note that $N=0$ iff $\ell_c=\ell$ and $N$ is even
iff $\ell_c(v)=1$. By Propositions \ref{qsenzill}, \ref{v=w} and
Theorem \ref{taulasn}, we have
$$\renewcommand\arraystretch{1.2}
\begin{array}{ll}
N=0\imp \sn(Q)=\pm,\ &\sn(\tilde{Q})=0,\\ N=1 \imp \sn(Q)=0, \
&\sn(\tilde{Q})=+,\\N=2 \imp \sn(Q)=0, \ &\sn(\tilde{Q})=0,\\
N=3 \imp \sn(Q)=0, \ &\sn(\tilde{Q})=-.
\end{array}\renewcommand\arraystretch{1.}
$$
There are 8 possibilities for $\ell_c$, one with $N=0$ or $N=3$
and three with $N=1$ or $N=2$. By Lemma \ref{nbrclas}, the
different values of $(c,d)$ lead to 2, 3, 3, 1 $k$-isomorphism
classes corresponding respectively to $N=0,1,2,3$. The number of
points of these curves is respectively $(N_1,N_2)=(q+1\pm
2\sqrt{2q},q^2+1),\,(q+1,q^2+1+4q),\,(q+1,q^2+1),\,(q+1,q^2+1-4q)$.
Thus, we get a contribution of respectively $(q-2)/6$, $(q-2)/6$,
$(q-2)/2$, $(q-2)/2$, $(q-2)/6$ curves in each of the sets
$\c{2\sqrt{2q}}{4q}$, $\c{-2\sqrt{2q}}{4q}$, $\c{0}{2q}$,
$\c{0}{0}$, $\c{0}{-2q}$.\e

From now on we deal with the case $m$ even.

\subsection{$\mathbf{P_{ab}(x)}$ irreducible}
By Proposition \ref{wp} we have $w=\tilde{w}=0$ and
$W=\tilde{W}=\{0\}$. We have thus $\sn(Q)=\pm$, and
$\sn(\tilde{Q})=+$, by Theorem \ref{taulasn}. On the other hand,
since $E_{ab}\colon k\lra k$ has a trivial kernel, we have
$E_{ab}(k)=k$ and we can always assume that $c=0$ by
(\ref{modiso}).

Let us study first the case $b=0$. By Proposition \ref{beso},
$P_{a0}(x)$ is irreducible iff $m\equiv 0\md4$ and $a\not\in
(k^*)^5$. In this case we have 8 $k$-isomorphism classes
represented by $(a,0,0,0)$, $(a,0,0,d_0)$, where $a$ runs on the 4
non-trivial classes of $k^*/(k^*)^5$. They have $(N_1,N_2)=(q+1\pm
\sqrt{q},q^2+1+q)$ and they contribute with 4 curves in each of
the sets $\c{-\sqrt{q}}{q}$, $\c{\sqrt{q}}{q}$.

In the case $a=b\ne0$ there are $\frac 25(q+1-[2]_{4\mid m})$
values of $a$ leading to $P_{aa}(x)$ irreducible, by Corollary
\ref{nbra}. As before, for each fixed value of $a$ we obtain 1
curve in each of the sets $\c{-\sqrt{q}}{q}$, $\c{\sqrt{q}}{q}$.

We have altogether a contribution of $\frac 25(q+1+[8]_{4\mid m})$
$k$-isomorphism classes in each of the sets $\c{-\sqrt{q}}{q}$,
$\c{\sqrt{q}}{q}$.

\subsection{$\mathbf{P_{ab}(x)=(1)(1)(3)}$}
By Proposition \ref{wp} we have $w=\tilde{w}=2$ and
$W=\gen{z_1,z_2}$, where $z_1,z_2$ are the roots of $P_{ab}(x)$ in
$k$. By Proposition \ref{beso}, we have necessarily $b\ne0$ and we
assume $a=b$. By Corollary \ref{nbra} we have $(q-1)/3$ values of
$a$ leading to this factorization of $P_{aa}(x)$. For any such
fixed value of $a$, the linear form $\ell$ on $W$ vanishes.

For any $c\in k$, let $N$ be the number of $z_i$ such that
$\ell_c(z_i)=0$. Note that $N=2$ iff $\ell_c=\ell$. By
Propositions \ref{qsenzill}, \ref{v=w} and Theorem \ref{taulasn},
we have
$$\renewcommand\arraystretch{1.2}
\begin{array}{ll}
N=0\imp \sn(Q)=0,\ &\sn(\tilde{Q})=-,\\ N=1 \imp \sn(Q)=0, \
&\sn(\tilde{Q})=+,\\N=2 \imp \sn(Q)=\pm, \ &\sn(\tilde{Q})=+.
\end{array}\renewcommand\arraystretch{1.}
$$There are 4
possibilities for $\ell_c$, one with $N=0$ or $N=2$ and two with
$N=1$. By Lemma \ref{nbrclas}, the different values of $(c,d)$
lead to 1, 2, 2 $k$-isomorphism classes corresponding respectively
to $N=0,1,2$. The number of points of these curves is respectively
$(N_1,N_2)=(q+1,q^2+1-2q),\,(q+1,q^2+1+2q),\,(q+1\pm
2\sqrt{q},q^2+1+2q)$. Thus, we get a contribution of respectively
$(q-1)/3$, $2(q-1)/3$, $(q-1)/3$, $(q-1)/3$ curves in each of the
sets $\c{0}{-q}$, $\c{0}{q}$, $\c{-2\sqrt{q}}{3q}$,
$\c{2\sqrt{q}}{3q}$.

\subsection{$\mathbf{P_{ab}(x)=(1)(2)(2)}$}
By Proposition \ref{wp} we have $w=2$, $\tilde{w}=4$ in this case.
If $P_{ab}(x)=(x+z)(x^2+v_1x+t_1)(x^2+v_2x+t_2)$, we have
$v_1+v_2=z$ and $W=\{0,z,v_1,v_2\}$.

Let us study first the case $b=0$. By Proposition \ref{beso},
$P_{a0}(x)=(1)(2)(2)$ iff $4\nmid m$. In this case, $(k^*)^5=k^*$
and we can assume $a=1$ by (\ref{modiso}). The linear form $\ell$
on $W$ is determined by $\ell(v_1)=\ell(v_2)=1$.

For any $c\in k$, let $N$ be the number of $v_i$ such that
$\ell_c(v_i)=0$. Note that $N=0$ iff $\ell_c=\ell$. By
Propositions \ref{qsenzill}, \ref{v=w} and Theorem \ref{taulasn},
we have
$$\renewcommand\arraystretch{1.2}
\begin{array}{ll}
N=0\imp \sn(Q)=\pm,\ &\sn(\tilde{Q})=0,\\ N=1 \imp \sn(Q)=0, \
&\sn(\tilde{Q})=0,\\N=2 \imp \sn(Q)=0, \ &\sn(\tilde{Q})=+.
\end{array}\renewcommand\arraystretch{1.}
$$There are 4
possibilities for $\ell_c$, one with $N=0$ or $N=2$ and two with
$N=1$. By Lemma \ref{nbrclas}, the different values of $(c,d)$
lead to 2,2,1 $k$-isomorphism classes according to $N=0,1,2$. The
number of points of these curves is respectively
$(N_1,N_2)=(q+1\pm
2\sqrt{q},q^2+1),\,(q+1,q^2+1),\,(q+1,q^2+1+4q).$
 Thus, we get respectively 1,1,2,1 curves in each of the sets $\c{-2\sqrt{q}}{2q}$,
 $\c{2\sqrt{q}}{2q}$, $\c{0}{0}$,
$\c{0}{2q}$.

In the case $a=b\ne0$ there are $(q/4)-[1]_{4\nmid m}$ values of
$a$ leading to $P_{aa}(x)=(1)(2)(2)$, by Corollary \ref{nbra}. As
before, for any such fixed value of $a$ we get respectively
1,1,2,1 curves with the same zeta functions as above.

We have altogether a contribution of $q/4$, $q/4$, $q/2$, $q/4$
$k$-isomorphism classes in each of the sets $\c{-2\sqrt{q}}{2q}$,
 $\c{2\sqrt{q}}{2q}$, $\c{0}{0}$,
$\c{0}{2q}$.

\subsection{$\mathbf{P_{ab}(x)=(1)(1)(1)(1)(1)}$}
By Proposition \ref{wp} we have $w=\tilde{w}=4$ in this case and
$W=\gen{z_1,z_2,z_3,z_4}$, where $z_1,z_2,z_3,z_4,z_1+z_2+z_3+z_4$
are the roots of $P_{ab}(x)$ in $k$.

Let us study first the case $b=0$. By Proposition \ref{beso},
$P_{a0}(x)$ splits completely in $k[x]$ iff $4\mid m$ and $a\in
(k^*)^5$. In this case we can assume $a=1$ by (\ref{modiso}), so
that $W=\mathbb{F}_{16}$. The linear form $\ell$ vanishes.

For any $c\in k$, let $N$ be the number of $z_i$ such that
$\ell_c(z_i)=0$. Note that $N=5$ iff $\ell_c=\ell$. By
Propositions \ref{qsenzill}, \ref{v=w} and Theorem \ref{taulasn},
we have
\begin{equation}\label{darrer}\renewcommand\arraystretch{1.2}
\begin{array}{ll}
N=1\imp \sn(Q)=0,\ &\sn(\tilde{Q})=-,\\ N=3 \imp \sn(Q)=0, \
&\sn(\tilde{Q})=+,\\N=5 \imp \sn(Q)=\pm, \ &\sn(\tilde{Q})=-.
\end{array}\renewcommand\arraystretch{1.}
\end{equation}
We cannot apply now Lemma \ref{nbrclas} but it is easy to
check that the different values of $(c,d)$ lead to 1,2,2
$k$-isomorphism classes according to $N=1,3,5$. The number of
points of these curves is respectively
$(N_1,N_2)=(q+1,q^2+1-4q),\,(q+1,q^2+1+4q),\,(q+1\pm
4\sqrt{q},q^2+1-4q).$
 Thus, we get respectively 1,2,1,1 curves in each of the sets $\c{0}{-2q}$,
$\c{0}{2q}$, $\c{-4\sqrt{q}}{6q}$, $\c{4\sqrt{q}}{6q}$.

In the case $a=b\ne0$ there are $\frac{q-4}{60}-[\frac 15]_{4\mid
m}$ values of $a$ leading to $P_{aa}(x)=(1)(1)(1)(1)(1)$, by
Corollary \ref{nbra}. For any such fixed value of $a$
(\ref{darrer}) holds. There are 16 possibilities for $\ell_c$,
five with $N=1$, ten with $N=3$ and one with $N=5$. By Lemma
\ref{nbrclas}, we get respectively 5,10,1,1 curves with the same
zeta functions as above.

We have altogether a contribution of $\frac{q-4}{12}$,
$\frac{q-4}6$, $\frac{q-4}{60}+[\frac 45]_{4\mid m}$,
$\frac{q-4}{60}+[\frac 45]_{4\mid m}$ $k$-isomorphism classes in
each of the sets $\c{0}{-2q}$, $\c{0}{2q}$, $\c{-4\sqrt{q}}{6q}$,
$\c{4\sqrt{q}}{6q}$.

\section{Jacobians in isogeny classes of supersingular abelian surfaces}
\label{isg}Let $A$ be a supersingular abelian surface defined over
$k$. Let $f_A(t)=t^4+a_1t^3+a_2t^2+qa_1t+q^2$ be the
characteristic polynomial of its Frobenius endomorphism. The
$k$-isogeny class of $A$ is determined by this polynomial, that
is, by the couple of integers $(a_1,a_2)$.

It is easy to list all couples $(a_1,a_2)$ that correspond to
supersingular abelian surfaces. The $k$-simple supersingular
isogeny classes can be found in \cite[Table 1]{mn}. If $A$ is
$k$-isogenous to a product of two supersingular elliptic curves,
we have $f_A(t)=f_{E_1}(t)f_{E_2}(t)$, with
$f_{E_i}(t)=t^2+b_it+q$. This gives,
$$
a_1=b_1+b_2,\quad a_2=2q+b_1b_2.
$$On the other hand, the possibilities for the integers $b_i$ were
determined by Waterhouse in his thesis \cite[Theorem 4.1]{wa}:
$$
b_i\in\{0,\,\pm \sqrt{2q}\},\,  \mbox{ if $m$ is odd};\quad\
b_i\in\{0,\,\pm \sqrt{q},\,\pm 2\sqrt{q}\},\, \mbox{ if $m$ is
even}.
$$
This gives respectively 6,15 $k$-split supersingular isogeny
classes according to $m$ being odd or even.

Gathering the computations of the previous section, we give in the
tables below the number of $k$-isomorphism classes of
supersingular curves of genus 2 whose jacobian lies in each
$k$-isogeny class. When in a column indexed as $(\pm a_1,a_2)$ we
say that $|\c{a_1}{a_2}|=N$, we mean that
$|\c{a_1}{a_2}|=|\c{-a_1}{a_2}|=N$.\be\e

\renewcommand\arraystretch{1.2}
\centerline{\begin{tabular}{|c|c|c|c|c|} \hline
$(a_1,a_2)$&$(0,0)$&$(0,2q)$&$(\pm\sqrt{2q},2q)$&$(\pm2\sqrt{2q},4q)$\\
\hline $|\c{a_1}{a_2}|$&$q-1$&$(q-2)/2$&$q/2$&$(q-2)/6$\\
\hline
$b_1,\,b_2$&$\sqrt{2q},\,-\sqrt{2q}$&$0,\,0$&$0,\,\pm\sqrt{2q}$&$b_1=b_2=\pm\sqrt{2q}$\\
\hline
\end{tabular}
}\e \centerline{\it Table 1: split isogeny classes, $m$ odd}\be\e

\centerline{\begin{tabular}{|c|c|c|c|c|} \hline
$(a_1,a_2)$&$(0,-2q)$&$(0,-q)$&$(0,q)$&$(\pm\sqrt{2q},q)$
\\\hline $|\c{a_1}{a_2}|$&$(q-2)/6$&0&$(q+1)/3$&$(q+1)/3$\\\hline
\end{tabular}
}\e \centerline{\it Table 2: simple isogeny classes, $m$ odd}\be\e

\centerline{\begin{tabular}{|c|c|c|c|c|c|} \hline
$(a_1,a_2)$&$(0,-2q)$&$(\pm\sqrt{q},0)$&$(0,q)$&$(0,2q)$&$(\pm\sqrt{q},2q)$\\
\hline $|\c{a_1}{a_2}|$&$(q-4)/12$&$0$&$2(q-1)/3$&$(5q-8)/12$&$0$\\
\hline
$b_1,\,b_2$&$2\sqrt{q},\,-2\sqrt{q}$&$\pm(\sqrt{q},\,-2\sqrt{q})$&$\sqrt{q},\,-\sqrt{q}$&$0,\,0$
&$0,\,\pm\sqrt{q}$\\
\hline
\end{tabular}
}\e \centerline{\begin{tabular}{|c|c|c|c|c|} \hline
$(a_1,a_2)$&$(\pm2\sqrt{q},2q)$&$(\pm2\sqrt{q},3q)$&$(\pm3\sqrt{q},4q)$&
$(\pm4\sqrt{q},6q)$\\
\hline $|\c{a_1}{a_2}|$&$q/4$&$(q-1)/3$&$0$&$\frac{q-4}{60}+[\frac45]_{4\mid m}$\\
\hline $b_1,\,b_2$&$0,\pm
2\sqrt{q}$&$b_1=b_2=\pm\sqrt{q}$&$\pm(\sqrt{q},\,2\sqrt{q})$
&$b_1=b_2=\pm2\sqrt{q}$\\
\hline
\end{tabular}
}\e \centerline{\it Table 3: split isogeny classes, $m$ even}\be\e

\centerline{\begin{tabular}{|c|c|c|c|} \hline
$(a_1,a_2)$&$(0,-q)$&$(0,0)$&$(\pm\sqrt{q},q)$
\\\hline $|\c{a_1}{a_2}|$&$(q-1)/3$&$q/2$&$\frac25(q+1+[8]_{4\mid m})$\\\hline
\end{tabular}
}\e \centerline{\it Table 4: simple isogeny classes, $m$
even}\be\e
\renewcommand\arraystretch{1.}

We see that some isogeny classes contain no jacobians. In most of
the cases there is a trivial explanation for this fact, but the
assertion that $\c{0}{-q}=\emptyset$ when $m$ is odd is far from
trivial and the achievement of this result was the initial
motivation for the paper.

\bigskip

\begin{tabular}{l}
 Daniel Maisner, Enric Nart\\
Departament de Matem\`atiques\\ Universitat Aut\`onoma de Barcelona\\
Edifici C\\ 08193 Bellaterra, Barcelona, Spain\\ {\it
danielm@mat.uab.es, nart@mat.uab.es}
\end{tabular}

\end{document}